
\input amstex
\documentstyle{amsppt}
\def\DJ{\leavevmode\setbox0=\hbox{D}\kern0pt\rlap
{\kern.04em\raise.188\ht0\hbox{-}}D}
\footline={\hss{\vbox to 2cm{\vfil\hbox{\rm\folio}}}\hss}
\nopagenumbers
\font\ff=cmr8
\def\txt#1{{\textstyle{#1}}}
\baselineskip=13pt
\def\hf{{\textstyle{1\over2}}}

\def\d{{\,\roman d}}
\def\e{\varepsilon}
\def\f{\varphi}
 \def\G{\Gamma}

\def\={\;=\;}

\def\zt{\zeta(\hf+it)}

\def\D{\Delta}

\def\z{\zeta}

\def\e{\varepsilon}
\def\D{\Delta}

\def\e{\varepsilon}

\font\teneufm=eufm10
\font\seveneufm=eufm7
\font\fiveeufm=eufm5
\newfam\eufmfam
\textfont\eufmfam=\teneufm
\scriptfont\eufmfam=\seveneufm
\scriptscriptfont\eufmfam=\fiveeufm
\def\mathfrak#1{{\fam\eufmfam\relax#1}}

\font\tenmsb=msbm10
\font\sevenmsb=msbm7
\font\fivemsb=msbm5
\newfam\msbfam
\textfont\msbfam=\tenmsb
\scriptfont\msbfam=\sevenmsb
\scriptscriptfont\msbfam=\fivemsb
\def\Bbb#1{{\fam\msbfam #1}}

\def \ZZ {\Bbb Z}

\def\rightheadline{{\hfil{\ff
Sums of squares of $|\zt|$}\hfil\tenrm\folio}}

\def\leftheadline{{\tenrm\folio\hfil{\ff
A. Ivi\'c }\hfil}}
\def\emptyheadline{\hfil}
\headline{\ifnum\pageno=1 \emptyheadline\else
\ifodd\pageno \rightheadline \else \leftheadline\fi\fi}

\topmatter
\title SUMS OF SQUARES OF $|\zt|$ OVER SHORT INTERVALS
\endtitle
\author   Aleksandar Ivi\'c \endauthor
\address{ \bigskip
Aleksandar Ivi\'c, Katedra Matematike RGF-a
Universiteta u Beogradu, \DJ u\v sina 7, 11000 Beograd,
Serbia (Yugoslavia).}
\endaddress
 \email {\tt aivic\@matf.bg.ac.yu,
aivic\@rgf.bg.ac.yu} \endemail \abstract {Sums of squares of
$|\zt|$ over short intervals are investigated. Known upper bounds
for the fourth and twelfth moment of $|\zt|$ are derived. A
discussion concerning other possibilities for the estimation of
higher power moments of $|\zt|$ is given.}
\endabstract
\endtopmatter

\noindent

\def\hf{{\textstyle{1\over2}}}
\def\txt#1{{\textstyle{#1}}}

\def\e{\varepsilon}
\def\f{\varphi}
\def\G{\Gamma}

\def\={\;=\;}

\def\zt{\zeta(\hf+it)}

\font\teneufm=eufm10
\font\seveneufm=eufm7
\font\fiveeufm=eufm5
\newfam\eufmfam
\textfont\eufmfam=\teneufm
\scriptfont\eufmfam=\seveneufm
\scriptscriptfont\eufmfam=\fiveeufm
\def\mathfrak#1{{\fam\eufmfam\relax#1}}

\font\tenmsb=msbm10
\font\sevenmsb=msbm7
\font\fivemsb=msbm5
\newfam\msbfam
      \textfont\msbfam=\tenmsb
      \scriptfont\msbfam=\sevenmsb
      \scriptscriptfont\msbfam=\fivemsb
\def\Bbb#1{{\fam\msbfam #1}}

\def \ZZ {\Bbb Z}

 \def\D{\Delta} \def\G{\Gamma} \def\e{\varepsilon}
 
\def\DJ{\leavevmode\setbox0=\hbox{D}\kern0pt\rlap
 {\kern.04em\raise.188\ht0\hbox{-}}D}

\heading 1. Introduction
\endheading
Let as usual $\z(s)$ denote the Riemann zeta-function.
The aim of this paper is to estimate the sum
$$
\sum_{r\le R}\int_{t_r-2G}^{t_r+2G}  \f_r(t)|\zt|^2\d t,\leqno(1.1)
$$
where $\f_r \in C^\infty$ is a non-negative function supported in
$[t_r-2G,\,t_r+2G]$ that equals unity in $[t_r-G,\,t_r+G]$,
$T \le t_1 < \cdots < t_R \le 2T$,
$t_{r+1} - t_r \ge 5G \,(r = 1,\ldots,R-1)$,
and $T^\e \le G \le T^{1/3}$ (here and later $\e$ denotes positive,
absolute constants, not necessarily the same ones at each occurrence).
This sum majorizes the classical sum
$$
\sum_{r\le R}\int_{t_r-G}^{t_r+G} |\zt|^2\d t,\leqno(1.2)
$$
which is of great importance in zeta-function theory (see K. Matsumoto
[18] for an extensive account on mean square theory involving $\z(s)$).
One can treat the sum in (1.1) by at least the following methods.

\hskip1cm
a) Using exponential averaging (or some other smoothing  like $\f_r$ above),
namely the Gaussian weight $\exp(-\hf x^2)$, in connection with ($\gamma$ is
Euler's constant)
$$
E(T) := \int_0^T|\zt|^2\d t - T\left(\log {T\over2\pi} + 2\gamma -1\right),
$$
since there exists a well-known  explicit formula of F.V. Atkinson [1]
for $E(T)$ (see also the author's monographs [5] and [6]).
This was done by D.R. Heath-Brown [2], who obtained
$$
\int_0^T|\zt|^{12}\d t \;\ll\;T^2\log^{17}T,\leqno(1.3)
$$
which is hitherto the best estimate of its kind, with many applications
to multiplicative number theory.

\hskip1cm
b) One can use the Voronoi summation formula (e.g., see [5, Chapter 3]) for
the explicit expression (approximate functional equation) for $|\zt|^2
= \chi^{-1}(\hf + it)\zeta^2(\hf + it)$, where
 $\z(s) = \chi(s)\zeta(1-s)$, namely
$$
\chi(s) \;=\;2^s\pi^{s-1}\sin(\hf\pi s)\G(1-s).
$$
Voronoi's formula is present indirectly in Atkinson's formula, so that
this approach is more direct. The effect of the smoothing function $\f_r$
in (1.1) is to shorten the sum approximating $|\z|^2$ to the range
${T\over2\pi}(1 - G^{-1}T^\e) \le n \le {T\over2\pi} \,(T = t_r)$.
After this no integration is needed, and proceeding as in [5, Chapters 7-8]
one obtains that the sum in (1.1) equals $O(RGT^\e)$ plus a multiple of
$$
\sum_{r\le R}\int_{t_r-2G}^{t_r+2G}  \f_r(t)
\sum_{k\le T^{1+\e}G^{-2}}(-1)^kd(k)k^{-1/2}\left({1\over4}
+ {t\over2\pi k}\right)^{-1/4}\sin f(t,k)\d t,\leqno(1.4)
$$
where $d(k)$ is the number of positive divisors of $k$ and
$$
f(t,k) := 2t\text{arsinh}\,\sqrt{\pi k\over2t} + \sqrt{2\pi kt + \pi^2k^2} -
\txt{1\over4}\pi, \text {arsinh}\, x = \log(x + \sqrt{1+x^2}).\leqno(1.5)
$$

\hskip1cm
c) Instead of the Voronoi summation formula one can use the (simpler)
Poisson summation formula, namely
$$
\sum_{n=1}^\infty f(n) = \int_0^\infty f(x)\d x +
2\sum_{n=1}^\infty  \int_0^\infty f(x)\cos(2\pi n x)\d x,\leqno(1.6)
$$
provided that $f(x)$ is smooth and compactly supported in $(0,\,\infty)$.
We shall give now a sketch of this approach. The integral in (1.1)
(with $t_r = T, \f_r(t) = \f(t)$) is majorized by $O(\log T)$ integrals
of the type
$$
I := \int_{T-2G}^{T+2G}\f(t)\Big|\sum_{N<n\le N'\le2N}n^{-1/2-it}\Big|^2\d t
\quad(G^{1+\e} \le N \ll \sqrt{T}, \f^{(j)}(t) \ll_j G^{-j}),
$$
since the contribution of $N \le G^{1+\e}$ is $\ll G^{1+\e}$ by the
mean value theorem for Dirichlet polynomials ([5, Theorem 5.2]). Squaring
out the above sum, estimating trivially the diagonal terms, and integrating
by parts sufficiently many times the contribution of the non-diagonal terms
we are led to the estimation of the expression
$$
\int_{T-2G}^{T+2G}\f(t)\sum_{1\le \ell\le N^{1+\e}G^{-1}}\sum
_{N<n\le N'\le2N}(n(n+\ell))^{-1/2}\left(1 + {\ell\over n}\right)^{it}\d t.
$$
Using Taylor's formula we can replace $(n+\ell)^{-1/2}$ by $n^{-1/2}$,
and then we write
$$
\sum_{N<n\le N'\le2N}n^{-1}\left(1 + {\ell\over n}\right)^{it}
= \sum_{N-G\le n\le N'+G}\Phi(n)n^{-1}\left(1 + {\ell\over n}\right)^{it}
+ O(GN^{-1}),\leqno(1.7)
$$
where $\Phi(x) \;(\ge 0)$ is a smooth function, supported in $[N-G,N_1+G]$,
equal to unity in $[N,N_1]$. This function facilitates truncation when
the Poisson summation formula (1.6) is applied; once with
$f(x) = \Phi(x)\cos\left(t\log\left(1 + {\ell\over x}\right)\right)$,
and once with the sine instead of cosine. The former gives the major
contribution equal to
$$\eqalign{&
\int_{T-2G}^{T+2G}\f(t)\sum_{\ell\le N^{1+\e}G^{-1}}
\sum_{m\le T^{1+\e}\ell N^{-2}}\times
\cr&
\times\int_{N-G}^{N_1+G}
\Phi(x)\cos\left(t\log\left(1 + {\ell\over x}\right)\right)\cos(2\pi mx)
\d x\,\d t.\cr}\leqno(1.8)
$$
The integral over $x$ is a linear combination of four exponential
integrals, one of which is
$$
\int_{N-G}^{N_1+G}\Phi(x)e^{iF(x)}\d x,\quad F(x) :=
t\log\left(1 + {\ell\over x}\right) + 2\pi mx,\leqno(1.9)
$$
so that
$$
F'(x) = t\left({1\over\ell+x} - {1\over x}\right) + 2\pi m,\quad
F''(x) = t\left({1\over x^2} - {1\over(\ell+x)^2}\right).
$$
The main contribution to (1.9) will come from the saddle point $x_0$ (see
e.g., [5, Chapter 2]) and will be a multiple of $|F''(x_0)|^{-1/2}
e^{iF(x_0)}$. This is the solution of $F'(x_0) = 0$, hence
$$
x_0 \;=\; {-2\pi m\ell + \sqrt{8\pi m\ell t + 4\pi^2m^2\ell}\over4\pi m}
\;\sim\; \sqrt{\ell t\over2\pi m}\qquad(T\to\infty),\leqno(1.10)
$$
where $x_0 \in [N,\,N_1]$, giving $m\asymp \ell TN^{-2}$, which is then
the summation condition in (1.8). Further we have
$F''(x_0) \sim t\ell x_0^{-3}$, hence the major contribution to (1.4) will
come from a multiple of
$$\eqalign{&
\int_{T-2G}^{T+2G}\f(t)t^{-1/4}\sum_{\ell\le N^{1+\e}G^{-1}}\ell^{-1/4}\cr&
\times\sum_{m\asymp\ell TN^{-2}}m^{-1/4}\exp\left(
it\left({\ell\over x_0} - {1\over2}{\ell^2\over x_0^2}
+ {1\over3}{\ell^3\over x_0^3}+ \ldots\right)+2\pi imx_0\right)\d t.\cr}
$$
Using (1.10) and grouping together the terms with $\ell m = k$, we
have $k \ll T^{1+\e}G^{-2}$, and after some transformations we arrive
essentially at the expression (1.4) with $R = 1, t_r = T$ and
$\left({1\over4}+ {t\over2\pi k}\right)^{-1/4}$ replaced by
$({t\over2\pi k})^{-1/4}$.

\bigskip
\heading
2. The twelfth moment of $|\zt|$
\endheading
The most direct way to treat (1.4) is first to develop $({1\over4} +
{t\over2k\pi})^{-1/4}$ by Taylor's formula. After this the mean value
theorem for integrals is applied. The essential part of (1.4) is then
$$
G\sum_{r\le R}\f_r(\tau_r)\tau_r^{-1/4}\sum_{k\le T^{1+\e}G^{-2}}
d(k)k^{-1/4}\sin f(\tau_r,k),\leqno(2.1)
$$
where $\tau_r \in [t_r-2G,\,\tau_r+2G]\,$, and (by considering
separately points with even and odd indices) we can assume the spacing
condition
$$
T/3 < \tau_1 < \cdots < \tau_R < 8T/3,\quad\tau_{r+1}-\tau_r \ge 5G\;\;
(r = 1,\ldots R-1).
$$
The important feature of (2.1) is that it contains no absolute
value signs, so that one can change the order of summation
and use the Cauchy-Schwarz inequality. Thus the expression in (2.1) equals
$$
\eqalign{&
G\sum_{k\le T^{1+\e}G^{-2}}d(k)k^{-1/4}
\sum_{r\le R}\f_r(\tau_r)\tau_r^{-1/4}\sin f(\tau_r,k)\cr&
\ll G\Bigl(\sum_{k\le T^{1+\e}G^{-2}}d^2(k)k^{-1/2}\Bigr)^{1/2}
\Bigl(\sum_{k\le T^{1+\e}G^{-2}}\Big|\sum_{r\le R}\f_r(\tau_r)
\tau_r^{-1/4}e^{if(\tau_r,k)}\Big|^2\Bigr)^{1/2}.
\cr}
$$
After squaring and changing the order of summation,
the last sum above becomes
$$
\eqalign{&
\sum_{r,s\le R}\f_r(\tau_r)\f_s(\tau_s)(\tau_r\tau_s)^{-1/4}
\sum_{k\le T^{1+\e}G^{-2}}\exp (if(\tau_r,k)-if(\tau_s,k))\cr&
\ll T^{-1/2}\sum_{r,s\le R}\Big|\sum_{k\le T^{1+\e}G^{-2}}
\exp (if(\tau_r,k)-if(\tau_s,k))\Big|\cr&
\ll_\e RT^{{1\over2}+\e}G^{-2} + T^{-1/2}\sum_{r,s\not= R}
\Big|\sum_{k\le T^{1+\e}G^{-2}}\exp (if(\tau_r,k)-if(\tau_s,k))\Big|.
\cr}
$$
The effect of this procedure is that there are no divisor coefficients in
the last sum over $k$. If the term $RT^{{1\over2}+\e}G^{-2}$ (coming from
the diagonal terms $r = s$) dominates, we would obtain not only (1.3),
but the stronger bound
$$
\int_0^T|\zt|^{6}\d t \;\ll_\e\;T^{1+\e},\leqno(2.2)
$$
which is not known to hold yet. For a conditional proof of (2.2), which
rests on a hypothesis involving the ternary additive divisor problem,
see the author's work [8]. However, the
estimation of the crucial exponential sum
$$
S \;:=\;\sum_{k\le T^{1+\e}G^{-2}}\exp (if(\tau_r,k)-if(\tau_s,k))
\leqno(2.3)
$$
is limited by the scope of the present-day exponential sum techniques
(see e.g., M.N. Huxley [4]),
and it does not seem that (2.2) can be reached (unconditionally)
in this fashion. Observing that
$$
{\partial f(x,k)\over\partial x} \,=\, 2\text{arsinh}\,\sqrt{\pi k
\over2x} \;\sim\; \sqrt{2\pi k\over x}\quad(x \asymp T),\leqno(2.4)
$$
setting
$$
f(u) := f(\tau_r,u) - f(\tau_s,u),\quad F := |\tau_r-\tau_s|
(KT)^{-1/2},
$$
we have (see [5, Chapters 1-2] for the relevant exponent pair technique)
that $S$ is split into $O(\log T)$ subsums of the type
$$\eqalign{&
\sum_{K<k\le K'\le2K}\exp(if(k)) \ll F^\kappa K^\lambda + F^{-1}\cr&
\ll J^{\kappa}T^{-\kappa/2}K^{\lambda-\kappa/2} + (KT)^{1/2}
|\tau_r-\tau_s|^{-1}\quad(K \le T^{1+\e}G^{-2}),\cr}
$$
provided that $|\tau_r-\tau_s| \le J (\ll T)$, and $(\kappa,\lambda)$ is
a (one-dimensional) exponent pair. Choosing $(\kappa,\lambda) =
(\hf,\hf)$, $J = T^{-\e}G^3$ we obtain (1.3) (with $T^\e$ in place of
$\log^{17}T$). Namely with $J = T^{-\e}G^3$ the number
of points $R = R_0$ to be estimated
satisfies $R_0 \ll_\e T^{1+\e}G^{-3}$, hence dividing $[T/2,\,5T/2]$
into subintervals of length not exceeding $J$ one obtains
$$
R \;\ll\; R_0(1 + T/J) \;\ll_\e\; T^{2+\e}G^{-6}
\;\ll_\e\; T^{2+\e}V^{-12},
$$
which easily yields (1.3) (with $T^\e$ in place of
$\log^{17}T$). This analysis was carried in detail in [5, Chapter 8],
where the possibilities of choosing other exponent pairs besides
$(\kappa, \lambda) = (\hf,\hf)$  were discussed. The method resembles
the method used originally by Heath-Brown [2] in his proof of (1.3).
\bigskip
\heading
3. The fourth moment of $|\zt|$
\endheading
We now present another approach. Let us start from the sum
$$
\sum := \sum_{r\le R}\int_{t_r-2G}^{t_r+2G}  \f_r(t)t^{-1/4}
\sum_{k\le T^{1+\e}G^{-2}}(-1)^kd(k)k^{-1/4}\sin f(t,k)\d t,\leqno(3.1)
$$
which is the dominant part of (1.4) after the removal of $({1\over4}
+ {t\over2\pi k})^{-1/4}$ by Taylor's formula. By H\"older's inequality
we obtain, if $M \in \Bbb N$ is fixed,
$$
\eqalign{
\sum &\ll \sum_{r\le R}T^{-1/4}\left(\int_{t_r-2G}^{t_r+2G}\f_r(t)
\Big|\sum_{k\le T^{1+\e}G^{-2}}(-1)^k\ldots\Big|^{2M}\d t\right)^{1\over2M}
G^{1-{1\over2M}}\cr&
\ll T^{-1/4}G^{1-{1\over2M}}\left(\sum_{r\le R}\int_{t_r-2G}^{t_r+2G}
\Big|\sum_{k\le T^{1+\e}G^{-2}}(-1)^k\ldots\Big|^{2M}\d t\right)^{1\over2M}
R^{1-{1\over2M}}\cr&
\ll T^{-1/4}(RG)^{1-{1\over2M}}\left(\int_{T/2}^{5T/2}\f(t)
\Big|\sum_{k\le T^{1+\e}G^{-2}}(-1)^k\ldots\Big|^{2M}\d t\right)^{1\over2M},
\cr}\leqno(3.2)
$$
since the intervals $[t_r-2G,\,t_r+2G]$ are non-overlapping. Here $\f(t)$
is a non-negative, smooth function supported in $[T/2,5T/2]$ such that
$\f(t) = 1$ for $T \le t \le 2T$, hence $\f^{(m)}(t) \ll_m T^{-m}$.

Take now $M = 1$ in (3.2) and split the sum over $k$ into $\ll \log T$
subsums over $K < k \le K' \le 2K$,\, where $K \le T^{1+\e}G^{-2}$. Then
$$
\eqalign{&
\int_{T/2}^{5T/2}\f(t)\Big|\sum_{K<k\le K'}(-1)^kd(k)k^{-1/4}
e^{if(t,k)}\Big|^2\d t\cr&
= O(TK^{1/2}\log^4K) \cr&
+ \sum_{K<k\not = m\le K'}
(-1)^{k+m}d(k)d(m)(km)^{-1/4}\int_{T/2}^{5T/2}\f(t)
e^{if(t,k)-if(t,m)}\d t\cr&
\ll_\e TK^{1/2}\log^4K + K^{\e-1/2}\sum_{K<k\not = m\le K'}
{K^{1/2}T^{1/2}\over|k-m|}\cr&
\ll_\e T^\e(TK^{1/2} + T^{1/2}K) \ll_\e T^{1+\e}K^{1/2},
\cr}
$$
where we used the first derivative test ([5, Lemma 2.1]) in conjunction with
(2.4). Thus for $M=1$ we obtain
$$
\sum \ll_\e T^{-1/4}(RG)^{1/2}T^{1/2+\e}K^{1/4} \ll_\e R^{1/2}T^{1/2+\e}.
\leqno(3.3)
$$
If $|\z(\hf + i\bar{t}_r)| \ge V > T^\e$ for $\bar{t}_{r+1} - \bar{t}_r
\ge 1 \;(r = 1,\ldots,R-1)$ and $\bar{t}_r\in [T,2T]$, then  on using
(see e.g., [6, Theorem 1.2])
$$
|\z(\hf + i\bar{t}_r)|^2 \ll \log T\left(\int_{\bar{t}_r-1}^{\bar{t}_r+1}
|\zt|^2\d t + 1\right)
$$
and grouping the points $\bar{t}_r$ into intervals of length $\ll G$
with center at $t_r$, choosing then $G = V^2T^{-\e}$, we obtain from (3.3)
$$
RV^2 \ll_\e R^{1/2}T^{1/2+\e},\quad R \ll_\e T^{1+\e}V^{-4},
$$
which easily gives then
$$
\int_0^T|\zt|^4\d t \;\ll_\e\; T^{1+\e}.\leqno(3.4)
$$
This is a weakened form of the fourth moment of $|\zt|$,
which is of the order $T\log^4T$. For an extensive account
on results on the fourth moment of $|\zt|$,
the reader is referred to [2], [7], [9]--[17], [19]--[21].
In particular, the sum (1.2) with $|\z|^4$ instead of $|\z|^2$   is
treated by spectral methods in [6], [15], [17] and [21]; the monograph
[21] of Y. Motohashi contains a comprehensive account of spectral theory.

\bigskip
\heading
4. Higher moments
\endheading

We start from (3.1), observing that we have, for
$t \asymp T,\,k \le T^{1+\e}G^{-2}$,
$$
f(t,k) = -\txt{1\over4}\pi + 2\sqrt{2\pi kt} +
\txt{1\over6}\sqrt{2\pi^3}k^{3/2}t^{-1/2} + a_5k^{5/2}t^{-3/2} +
a_7k^{7/2}t^{-5/2} + \ldots
\leqno(4.1)
$$
with effectively computable constants $a_{2\ell-1}\;(\ell \ge 3)$. Noting
that for $G = V^2T^{-\e}$
$$
k^{5/2}t^{-3/2} \;\ll\;T^{1+\e}G^{-5} \;\le\;T^{-\e}
$$
for
$$
V \;\ge\;T^{{1\over10}+\e},\leqno(4.2)
$$
it follows that we shall obtain (provided that (4.2) holds) (3.1) with
$f(t,k)$ replaced by
$$
-\txt{1\over4}\pi + 2\sqrt{2\pi kt} +
\txt{1\over6}\sqrt{2\pi^3}k^{3/2}t^{-1/2}
$$
times a series whose terms are of descending order of magnitude.
The main contribution will thus come from the above term, and proceeding
as in the proof of (3.4) we obtain

THEOREM 1. {\it Let $T \le t_1 < \ldots < t_R \le 2T$ be points
such that $|\z(\hf + it_r)| \ge VT^{-\e}$ with $t_{r+1} - t_r \ge
V \ge T^{{1\over10}+\e}$ for $r = 1,\ldots, R-1$. Then, for any
fixed integer $M\ge 1$,} $$\eqalign{& R \ll_\e
T^{\e-M/2}V^{-2}\max_{K\le T^{1+\e}V^{-4}}\times\cr&\times
\int_{T/2}^{5T/2}\f(t)\Big|\sum_{K\le k\le
K'\le2K}(-1)^kd(k)k^{-1/4} \exp(2i\sqrt{2\pi
kt}+cik^{3/2}t^{-1/2})\Big|^{2M}\d t,\cr}\leqno(4.3) $$ {\it where
$c = \sqrt{2\pi^3}/6$ and $\f(t)$ is a non-negative, smooth
function supported in \break $[T/2,5T/2]$ such that $\f(t) = 1$
for $T\le t \le2T$}.

\bigskip
In Section 3 we investigated the case $M=1$ of (4.3),
and it was shown that this leads
to the known bound (3.4). The next case $M=2$ is probably the most promising
one, and perhaps could lead to a sharper bound than (1.3). It reduces to the
evaluation of the integral
$$\eqalign{&
\int_{T/2}^{5T/2}\f(t)\sum_{K<m,n,k,\ell\le K'}(-1)^{m+n+k+\ell}\times\cr&
\times d(m)d(n)d(k)d(\ell)(mnk\ell)^{-1/4}\exp(iDt^{1/2}+iEt^{-1/2})\d t,\cr}
\leqno(4.4)
$$
where we have set
$$\eqalign{
&D \= D(m,n,k,l) :\= 2\sqrt{2\pi}(\sqrt{m}+\sqrt{n}-\sqrt{k}-\sqrt{\ell}),
\cr&
E \= E(m,n,k,l) :\= \txt{1\over6}\sqrt{2\pi^3}
(m\sqrt{m}+n\sqrt{n}-k\sqrt{k}-\ell\sqrt{\ell}).\cr}
$$
This problem bears resemblance to the evaluation of
$$
\int_1^X\D^4(x)\d x,\;\D(x) \= \sum_{n\le x}d(n) - x(\log x + 2\gamma-1),
$$
where $\D(x)$ represents the error term in the classical Dirichlet
divisor problem. This was investigated by K.-M. Tsang [23], and the
connection between the fourth moments of $\D(x)$ and
$E(T)$ was considered by the author
[12]. On one hand the problem considered by Tsang was more difficult,
since he aimed at an asymptotic formula, and not only an upper bound.
On the other hand, the presence of $iEt^{-1/2}$ in the exponential
in (4.4) induces extra difficulties.

If $D=0$ in (4.4), then the integral is estimated trivially as $O(T)$.
One has (see [22]) $D=0$ if $(m,n) = (k,\ell)$ or
$(m,n) = (\ell,k)$ or $m = \alpha^2h, n = \beta^2h, k = \gamma^2h,
\ell = \delta^2h, \alpha + \beta = \gamma + \delta$ and $h$ is
squarefree. Thus the contribution of $(m,n,k,\ell)$ for which $D=0$ is
$$
\eqalign{&
\ll T\Bigl(\bigl(K^{-1}\sum_{K<k\le K'}d^2(k)\bigr)^2 +
K^{\e-1}\sum_{h=1}^\infty\,\sum_{\alpha,\beta,\gamma\asymp(K/h)^{1/2}}1
\Bigr)\cr&
\ll_\e\;T^{1+\e}K.\cr}
$$
When inserted in (4.4) this portion will yield ((4.3) with $M=2$)
$$
R \ll_\e T^{\e-1}V^{-2}\cdot T^{1+\e}\max_{K\le T^{1+\e}V^{-4}}K
\ll_\e T^{1+\e}V^{-6},\leqno(4.5)
$$
and (4.5) is the estimate that leads to (2.2). As in the discussion
involving the twelfth moment estimate of  $|\zt|$, it turns out that
non-diagonal terms, i.e. those for which $D \not = 0$ in our case, are
the difficult ones to estimate. We suppose now that $D>0$ (the case
$D<0$ is analogous) and $E>0$. Namely in case $E\le 0$ the estimation
is easier. When $E < 0$ the relevant exponential integral has no saddle
point and the contribution of $D > T^{\e-1/2}$ is negligible.
If $E=0$ the estimation is much easier. Write now (4.4) as
$$
\sum_{K<m,n,k,\ell\le K'}(-1)^{m+n+k+\ell}d(m)d(n)d(k)d(\ell)
(mnk\ell)^{-1/4}I(T),
$$
where
$$
I(T) = I(T;m,n,k,\ell) := \int_{T/2}^{5T/2}\f(t)
e^{iDt^{1/2} + iEt^{-1/2}}\d t.\quad
$$
Change of variable ($x = \sqrt{t},\,\Phi(x) = 2x\f(x^2)$) gives
$$
I(T) = \int^{\sqrt{5T/2}}_{\sqrt{T/2}}\Phi(x)e^{iDx+iE/x}\d x,
\quad \Phi^{(j)}(x) \ll_j x^{1-j}\;(j = 0,1,2,\ldots),\leqno(4.6)
$$
and it is assumed that $D,E > 0$. If
$$
D \le \eta K^{-1/2} \leqno(4.7)
$$
for a fixed, small $\eta > 0$, then from
$\sqrt{\ell} = \sqrt{m}+\sqrt{n}-\sqrt{k} + O(D)$ it follows by squaring that
$$
\ell = m + n + k + 2(\sqrt{mn} - \sqrt{mk} - \sqrt{nk}) + O(DK^{1/2}).
\leqno(4.8)
$$
Hence for suitable small $\eta$ the $O$-term in (4.8) will be $< 1/3$
in absolute value, so that there will be at most one $\ell$ for any given
choice of $m,n,k$. Moreover, for fixed $m,n$ it follows from (4.8) that
we must have ($||x||$ is the distance of $x$ to the nearest integer)
$$
||\,2\sqrt{k}(\sqrt{m} + \sqrt{n}) - 2\sqrt{mn}\,|| \ll DK^{1/2}\leqno(4.9)
$$
for any $k$ for which there will exist an $\ell$ satisfying (4.8).
If $D > \eta K^{-1/2}$, then we use the first derivative test and
trivial estimation to obtain a bound sharper then (4.15)

Suppose that beside (4.7) one also has
$$
D > T^{\e-1/2}.\leqno(4.10)
$$
From
$$\eqalign{
I(T) &= i\int^{\sqrt{5T/2}}_{\sqrt{T/2}}
{\left({\Phi(x)\over D - Ex^{-2}}\right)}'e^{iDx+iE/x}\d x\cr&
= i\int^{\sqrt{5T/2}}_{\sqrt{T/2}}
\left(\Phi'(x) - {2Ex^{-3}\over D - Ex^{-2}}\Phi(x)\right){e^{iDx+iE/x}
\over D - Ex^{-2}}\,\d x
\cr}
$$
we see that, if $D > C_2E/T$ or $D < C_1E/T$ for suitable constants
$0 < C_1 < C_2$, then repeated integrations by parts show, in view
of (4.10), that the contribution of $I(T)$ is negligible.
If $C_1E \le DT \le C_2E$, then by the second derivative
test ([5, Lemma 2.2]) we have
$$
I(T) \;\ll\; E^{3/4}D^{-5/4},\leqno(4.11)
$$
and this bound is non-trivial if (4.10) holds.
The number of $k$'s satisfying
$$
||\,2\sqrt{k}(\sqrt{m} + \sqrt{n}) - 2\sqrt{mn}\,|| \;<\;\delta
\quad(0<\delta < \hf)\leqno(4.12)
$$
for a fixed $m,n$ is uniformly
$$
\ll\;K\delta + K^{2/3}\leqno(4.13)
$$
by [22, Lemma 4]. This bound follows by applying standard exponential
sum techniques. Namely $||x|| < \delta\;(0 < \delta < 1/2)$ is
equivalent to
$$
[x+\delta] - [x-\delta] = \psi(x-\delta) - \psi(x+\delta) + 2\delta = 1
$$
(the expression equals zero otherwise) with
$$
\eqalign{\psi(x) &= x - [x] - \hf = -{1\over\pi}\sum_{n=1}^\infty
{\sin(2\pi nx)\over n}\cr& = -{1\over\pi}\sum_{n=1}^N
{\sin(2\pi nx)\over n}
+ O\left(\min\left(1,\,{1\over N||x||}\right)\right)
\quad(x\not\in \ZZ,\, N \ge 3).\cr}
$$
Thus the problem is quickly reduced to a problem involving the
estimation of exponential sums.
In our case $D \ll E/T \ll K^{3/2}/T$, hence (in
view  of (4.9)) in (4.12) we may
take $\delta = CDK^{1/2} \ll K^2/T$ with a suitable $C>0$.
The contribution for which (4.7) and (4.10) hold is, on using (4.11) and
(4.13),
$$\eqalign{
R &\ll_\e T^{\e-1}G^{-1}E^{3/4}D^{-5/4}K(K^{3/2}D + K^{2/3})\cr&
\ll_\e T^\e(T^{-3/4}K^{13/4}G^{-1} + T^{-1}G^{-1}(E/D)^{3/4}K^{5/3}D^{-1/2})
\cr&
\ll_\e T^\e(T^{-3/4}(TG^{-2})^{13/4}G^{-1} + T^{-1}G^{-1}T^{3/4}(TG^{-2})^{5/3}
T^{1/4})\cr&
\ll_\e T^{{5\over2}+\e}V^{-15} + T^{{5\over3}+\e}V^{-{26\over3}}.
\cr}                        \leqno(4.14)
$$

If $D > \eta K^{-1/2}$, then (4.10) holds. From (4.8) it follows that there
are $\ll DK^{1/2}$ possible values of $\ell$, hence by trivial estimation,
similar to the first bound in (4.14), we obtain a contribution which is
$\ll_\e T^{5/2+\e}V^{-15}$.

If $D < T^{\e-1/2}$, then the bound in (4.13) and trivial estimation give
$$\eqalign{
R &\ll_\e T^{\e-1}G^{-1}KT(K^{3/2}T^{-1/2} + K^{2/3})\cr&
\ll_\e T^{2+\e}V^{-12} + T^{{5\over3}+\e}V^{-{26\over3}}.\cr}\leqno(4.15)
$$
However, the term $K^{2/3}$ in (4.13) can be replaced, by applying
e.g. the Bombieri-Iwaniec method (see e.g., M.N. Huxley [4])
for the estimation of the relevant
exponential sum, by the sharper term $K^{7/11}$, yielding
$$
R \ll_\e T^{2+\e}V^{-12} + T^{{18\over11}+\e}V^{-{94\over11}}.
$$
Noting that $V \le T^{1/6}$ implies that $T^2V^{-12} \le T^{5/2}V^{-15}$,
we obtain (under (4.2))
$$
R \ll_\e T^{{5\over2}+\e}V^{-15} + T^{{18\over11}+\e}V^{-{94\over11}},
\leqno(4.16)
$$
which is a little weaker than the twelfth moment estimate
$$
R \ll_\e T^{2+\e}V^{-12}.\leqno(4.17)
$$
The above method does give $R \ll_\e T^{{5\over2}+\e}V^{-15}$, if
one uses the recent result of O. Robert--P. Sargos [22] (I am grateful
to P. Sargos for pointing this out to me).
Other improvements are also possible. Namely we used the bound in (4.13)
 assuming that $m,n$
are fixed. However, summation over $m,n$ in the relevant exponential
sum leading to (4.13)  produces in fact a three-dimensional exponential
sum which can be treated e.g. by three-dimensional exponent pairs.
In this way presumably one obtains $R \ll_\e T^{{5\over2}+\e}V^{-15}$
if (4.2) holds, however if (4.2) does not hold, then this bound easily
follows from (3.4). A further possibility for improvement is to
note that the saddle point method, in case $C_1E \le DT \le C_2E$,
furnishes in fact an asymptotic formula for $I(T)$, and not just the
upper bound in (4.11). This will lead to the exponential sum
$$ \eqalign{&
{\sum_{m,n,k,\ell}}^*(-1)^{m+n+k+\ell}
d(m)d(n)d(k)d(\ell)(mnk\ell)^{-1/4}\cr&\times\Phi\left(\sqrt{E\over D}
\right)E^{1/4}D^{-3/4}\exp(2i\sqrt{DE}),\cr}
\leqno(4.18)
$$
where ${}^*$ denotes that (4.9) and $E \asymp DT$ hold.
This sum is difficult to handle, since the variables in the exponential
are not separated and, for fixed $m,n,k$, there will exist at most one $\ell$.
 However, in view of $C_1E \le DT \le C_2E$, we
expect heuristically that certain separation of variables in (4.18)
is possible, which would lead to considerable cancellation. Presumably
this will yield at least (4.17).

\vfill\break

\Refs

\item{[1]} F.V. Atkinson, `The mean value of the zeta-function on
the critical line', {\it Quart. J. Math. Oxford} {\bf 10}(1939), 122-128.

\item{[2]}  D.R. Heath-Brown, `The twelfth power moment of the
Riemann zeta-function', {\it Quart. J. Math. Oxford} {\bf29}(1978), 443-462.

\item{[3]} { D.R. Heath-Brown}, `The fourth moment of the Riemann
zeta-function',
{\it Proc. London Math. Soc.}  (3){\bf38}(1979), 385-422.

\item{[4]} M.N. Huxley, `Area, lattice points, and exponential sums',
{\it LMS Monographs, New Ser. {\bf13}
Oxford University Press}, Oxford, 1996.

\item {[5]} { A. Ivi\'c},  `The Riemann zeta-function', {\it John Wiley
and Sons}, New York, 1985.

\item {[6]} { A. Ivi\'c},  `Mean values of the Riemann zeta-function',
LN's {\bf82}, {\it Tata Institute of Fundamental Research}, Bombay, 1991
(distr. by Springer Verlag, Berlin etc.).

\item{ [7]} { A. Ivi\'c},  `On the fourth moment of the Riemann
zeta-function', {\it Publs. Inst. Math. (Belgrade)}
{\bf57(71)}(1995), 101-110.

\item{[8]} { A. Ivi\'c},
`On the ternary additive divisor problem and the sixth moment of the
zeta-function', {\it ``Sieve Methods, Exponential Sums, and
their Applications in
Number Theory" }, Cambridge
University Press, Cambridge, 1996, 205-243.

\item{[9]} A. Ivi\'c, `On the error term for the fourth moment of the
Riemann zeta-function', {\it J. London Math. Soc.},
{\bf60}(2)(1999), 21-32.

\item{[10]} A. Ivi\'c,  `On the integral of the error term in
the fourth moment of the Riemann zeta-function', {\it Functiones
et Approximatio}, {\bf28}(2000), 37-48.

\item{[11]} A. Ivi\'c, `The Laplace transform of the fourth moment
of the zeta-function', {\it Univ. Beog. Publik. Elektroteh. Fak. Ser.
Matematika} {\bf11}(2000), 41-48.

\item{[12]} A. Ivi\'c, `On some problems involving the mean square
of $\zt$', {\it Bulletin CXXI Acad. Serbe Sci. Arts, Classe Sci. Math.}
{\bf 23}(1998), 71-76.

\item{ [13]}  A. Ivi\'c and Y. Motohashi,  `A note on the mean value of
the zeta and L-functions VII', {\it  Proc. Japan Acad. Ser. A}
{\bf  66}(1990), 150-152.

\item{ [14]}{ A. Ivi\'c and Y. Motohashi}, `The mean square of the
error term for the fourth moment of the zeta-function',
{\it Proc. London Math.
Soc.} (3){\bf66}(1994), 309-329.

\item {[15]} { A. Ivi\'c and Y. Motohashi},  `The fourth moment of the
Riemann zeta-function', {\it J. Number Theory} {\bf 51}(1995), 16-45.

\item {[16]}  A. Ivi\'c, M. Jutila and Y. Motohashi,
{\it The Mellin transform of powers of the zeta-function}, {\it Acta
Arith.} {\bf95}(2000), 305-342.

\item{ [17]}  H. Iwaniec, `Fourier coefficients of cusp forms
and the  Riemann zeta-function ', {\it S\'em. Th\'eorie Nombres Bordeaux},
Exp. No. {\bf18}, 1979/80.

\item{[18]} K. Matsumoto, `Recent developments in the mean square theory
of the Riemann zeta and other zeta-functions', in ``{\it Number Theory}",
Birkh\"auser, Basel, 2000, 241-286.

\item{ [19]} Y. Motohashi,   `An explicit formula for the fourth power
mean of the Riemann zeta-function', {\it Acta Math. }{\bf 170}(1993),
181-220.

\item {[20]} Y. Motohashi,  `A relation  between the Riemann
zeta-function and the hyperbolic Laplacian', {\it Annali Scuola Norm.
Sup. Pisa, Cl. Sci. IV ser.} {\bf 22}(1995), 299-313.

\item {[21]} Y. Motohashi,  `Spectral theory of the Riemann
zeta-function', {\it Cambridge University Press}, Cambridge, 1997.

\item{[22]} O. Robert and P. Sargos, `Three-dimensional
exponential sums with monomials', J. reine angew. Math. (in print).

\item {[23]} K.-M. Tsang, `Higher-power moments of $\D(x)$, $E(t)$ and
$P(x)$', {\it Proc. London Math. Soc.} (3){\bf65}(1992),   65-84.

\bigskip

Aleksandar Ivi\'c

 \par Katedra Matematike RGF-a

 Universiteta u
Beogradu\par \DJ u\v sina 7, 11000 Beograd

Serbia (Yugoslavia)\par  e-mail: {\tt aivic\@matf.bg.ac.yu,
aivic\@rgf.bg.ac.yu}

\endRefs

\bye